\documentclass{amsart}

\usepackage{fullpage}
\usepackage{amsfonts}
\usepackage{amssymb}

\input xy
\xyoption{matrix}\xyoption{arrow}\xyoption{curve}

\begin{document}

\newcommand{\aut}{\ensuremath{\mathrm{Aut}}}
\newcommand{\ind}{\ensuremath{\mathrm{ind}}}
\newcommand{\Hom}[3]{\ensuremath{\mathrm{Hom}_{#1}(#2,#3)}}
\newcommand{\stHom}[3]{\ensuremath{\underline{\mathrm{Hom}}_{#1}(#2,#3)}}
\newcommand{\Ext}[4]{\ensuremath{\mathrm{Ext}^{#2}_{#1}(#3,#4)}}
\newcommand{\Tor}[4]{\ensuremath{\mathrm{Tor}^{#1}_{#2}(#3,#4)}}
\newcommand{\coker}{\ensuremath{\mathrm{coker}}}
\newcommand{\fmod}[2][]{\ensuremath{\mathrm{mod}_{#1}\mbox{-} #2}}
\newcommand{\Mod}[2][]{\ensuremath{\mathrm{Mod}_{#1}\mbox{-} #2}}
\newcommand{\Bimod}[2][]{\ensuremath{\mathrm{Bimod}_{#1}\mbox{-} #2}}
\newcommand{\bimod}[2][]{\ensuremath{\mathrm{bimod}_{#1}\mbox{-} #2}}
\newcommand{\lmod}[1]{\ensuremath{#1 \mbox{-} \mathrm{mod}}}
\newcommand{\stmod}[1]{\ensuremath{\mathrm{\underline{mod}}\mbox{-} #1}}
\newcommand{\lstmod}[1]{\ensuremath{#1 \mbox{-} \mathrm{\underline{mod}} }}
\newcommand{\proj}[1]{\ensuremath{\mathrm{proj} \mbox{-} #1}}
\newcommand{\Proj}[1]{\ensuremath{\mathrm{Proj}\ #1}}
\newcommand{\inj}[1]{\ensuremath{\mathrm{inj} \mbox{-} #1}}
\newcommand{\modp}[1]{\ensuremath{\mathrm{mod}_{\mathcal{P}} #1}}
\newcommand{\ses}[5]{\ensuremath{0 \rightarrow #1 
\stackrel{#4}{\longrightarrow} 
#2 \stackrel{#5}{\longrightarrow} #3 \rightarrow 0}}
\newcommand{\add}[1]{\ensuremath{\mathrm{add}(#1)}}
\newcommand{\rad}{\ensuremath{\mathrm{rad}\ }}
\newcommand{\soc}{\ensuremath{\mathrm{soc}\ }}
\newcommand{\ann}{\ensuremath{\mathrm{ann}}}
\newcommand{\bimo}[1]{\ensuremath{{}_{#1}#1_{#1}}}
\newcommand{\pdim}{\ensuremath{\mathrm{pd}}}
\newcommand{\End}[2]{\ensuremath{\mathrm{End}_{#1}(#2)}}
\newcommand{\stEnd}[2]{\ensuremath{\underline{\mathrm{End}}_{#1}(#2)}}
\newcommand{\Tr}{\ensuremath{\mathrm{Tr}}}
\newcommand{\im}{\ensuremath{\mathrm{im}}}
\newcommand{\op}{\ensuremath{\mathrm{op}}}
\newcommand{\und}[1]{\underline{#1}}
\newcommand{\gen}[1]{\ensuremath{\langle #1 \rangle}}
\newcommand{\chr}{\ensuremath{\mathrm{char}}}
\newcommand{\HH}{\mathrm{HH}}

\newtheorem{therm}{Theorem}[section]
\newtheorem{defin}[therm]{Definition}
\newtheorem{propos}[therm]{Proposition}
\newtheorem{lemma}[therm]{Lemma}
\newtheorem{coro}[therm]{Corollary}

\title{Periodic resolutions and self-injective algebras of finite type}
\author{Alex S. Dugas}
\address{Department of Mathematics, University of California, Santa Barbara, CA 93106, USA}
\email{asdugas@math.ucsb.edu}


\begin{abstract}  We say that an algebra $A$ is periodic if it has a periodic projective resolution as an $(A,A)$-bimodule.  We show that any self-injective algebra of finite representation type is periodic.  To prove this, we first apply the theory of smash products to show that for a finite Galois covering $B \rightarrow A$, $B$ is periodic if and only if $A$ is.  In addition, when $A$ has finite representation type, we build upon results of Buchweitz to show that periodicity passes between $A$ and its stable Auslander algebra.  Finally, we use Asashiba's classification of the derived equivalence classes of self-injective algebras of finite type to compute bounds for the periods of these algebras, and give an application to stable Calabi-Yau dimensions.
\end{abstract}

\maketitle
 
\section{Introduction}

One of the aims of this article is to investigate periodicity of Hochschild cohomology for a finite-dimensional algebra $A$ over an algebraically closed field $k$.  Such periodicity is clearly guaranteed if the minimal projective resolution of the bimodule ${}_{A} A_{A}$ over the enveloping algebra $A^e = A^{\op} \otimes_k A$ is periodic, and we thus say that an algebra with this property is {\it periodic}.  It is shown in \cite{GrSnSo} that periodic algebras are necessarily self-injective, and numerous examples are known.  Schofield has shown that the preprojective algebras associated to Dynkin graphs are periodic (see \cite{ErdSna2}), and these results have recently been generalized to deformed preprojective algebras by Bia\l kowski, Erdmann and Skowro\'{n}ski \cite{DPA}.  By direct calculation, Erdmann, Holm and Snashall have also verified that the self-injective algebras of finite type and tree class $\mathbb{A}_n$ are periodic \cite{EH, EHS}.  Additional examples, including trivial extensions of path algebras of Dynkin quivers, were discovered by Brenner, Butler and King \cite{PAAK}.   A fairly comprehensive survey of these algebras with still more examples is given by Erdmann and Skowro\'nski in \cite{ErdSko2}.

Self-injective algebras of finite representation type provide a particularly interesting problem in this context.  It is easy to see that every nonprojective indecomposable module $M$ over such an algebra $A$ must be isomorphic to one of its syzygies.  Using this observation, Green, Snashall and Solberg show that some syzygy of $A$ over $A^e$ is isomorphic to a twisted bimodule ${}_1 A_{\sigma}$ for some $\sigma \in \aut(A)$ \cite{GrSnSo}.  This implies that the minimal projective resolution of $A$ over $A^e$ is very close to being periodic (for instance, the modules in the resolution repeat).  Nevertheless, whether or not $A$ is actually periodic has serious implications for the structure of the Hochschild cohomology ring of $A$.  In particular, if $\mathcal{N}$ denotes the nil radical of the Hochschild cohomology ring $\HH^*(A)$, Green, Snashall and Solberg show that $$\HH^*(A)/ \mathcal{N} \cong \left\{ \begin{array}{cc} k, & \mathrm{if}\ A\ \mathrm{is\ not\ periodic} \\ k[x], & \mathrm{if}\ A \ \mathrm{is\ periodic,} \end{array} \right.$$
where the degree of $x$ equals the (minimal) period of $A$.

In this article, we resolve this question, showing that self-injective algebras of finite type are indeed periodic.  Erdmann and Skowro\`nski have recently obtained this result for standard algebras not of type $(D_{3m}, s/3, 1)$ with $3 \nmid s$ by different means in \cite{ErdSko2}.  In the course of our proof, we expand the known ways of finding periodic algebras by showing that this periodicity is preserved upon passing between an algebra and its finite Galois coverings (equivalently, smash products), as well as its stable Auslander algebra when it has finite type.  In particular, these methods allow us to establish a connection between the periodicity of preprojective algebras and self-injective algebras of finite type.  Using this correspondence, we can even compute the periods of many algebras in the latter class in terms of their {\it types} as defined in \cite{DECSA}, and obtain decent bounds for the periods of the rest.  

Recently there has been renewed interest in periodicity questions arising from the study of Calabi-Yau dimensions of stable module categories \cite{BiaSko, ErdSko, CYFA}, and we apply our results to calculate these dimensions for the standard symmetric algebras of finite type.  This work corrects an error in \cite{ErdSko}, and suggests shortcomings of the proofs of \cite{BiaSko} concerning which finite-type self-injective algebras are Calabi-Yau.  We shall return to this problem in a subsequent paper.


I would like to extend thanks to Hideto Asashiba for some helpful correspondence, and to Kevin Walker whose examples motivated this work. 

\section{Preliminaries}

All algebras we consider are assumed to be split, basic finite-dimensional algebras over a field $k$. Such algebras can always be expressed as path algebras modulo relations $kQ/I$ for a quiver $Q$, and we will usually assume that we are given such a presentation.  In this case, we write $e_1, \ldots, e_n$ for the primitive idempotents associated to the vertices $Q_0$ of $Q$, and we write $Q_1$ for the set of arrows of $Q$.  We write $\fmod{A}$ (resp. $\Mod{A}$) for the category of finite-dimensional (resp. all) right $A$-modules, which we identify with contravariant representations of $Q$, and we denote the simple right $A$-modules (up to isomorphism) as $S_i$ for $1 \leq i \leq n$.  We let $A^e = A^{\op} \otimes_k A$ be the enveloping algebra for $A$, and we identify $(A,A)$-bimodules with right $A^e$-modules.

Suppose that $A = \oplus_{g \in G} A_g$ is a $G$-graded algebra for some group $G$ (with identity $e$).  We assume throughout this article, that the primitive idempotents $e_i$ of $A$ are homogeneous of degree $e$.  This ensures that the indecomposable projective modules are graded.  We write $J = J(A)$ for the Jacobson radical of $A$ and $J_G = J_G(A)$ for the graded Jacobson radical of $A$, which can be defined as the intersection of the maximal graded right ideals of $A$.  We will often assume that the grading is such that $J_G = J$, which we term a {\it radical grading}.  When $G$ is finite, Cohen and Montgomery have shown that $J_G$ is the largest homogeneous ideal contained in $J$ \cite{CoMo}; hence in this case, equality of the two radicals is tantamount to homogeneity of $J$.  Furthermore, this equality is automatic whenever $|G|$ is invertible in $A$.  One significant way of obtaining radical gradings of $A = kQ/I$ is through a function $\pi : Q_1 \rightarrow G$ as in \cite{Green}, whose image generates $G$ and for which $I \subset kQ$ is homogeneous.  Here $J$ is a homogeneous ideal containing $\oplus_{g \neq e} A_g$.  As illustrated in \cite{Green}, these gradings correspond to Galois covers of $A$.  

We write $\fmod[G]{A}$ for the category of finite-dimensional $G$-graded right $A$-modules and degree-preserving morphisms.  For a graded $A$-module $M$ and $d \in G$, we define $M[d]$ to be the graded $A$-module given by $M[d]_g = M_{d^{-1}g}$.  For a radical grading, the graded simple $A$-modules, up to isomorphism, are precisely the $S_i[d]$ for $1 \leq i \leq n$ and $d \in G$, and each such has a minimal graded projective resolution.  Furthermore, the minimal graded projective resolution of $S_i[d]$ coincides with the minimal projective resolution of $S_i$ as a complex of ungraded modules.

We shall also consider graded $(A,A)$-bimodules.  We say a bimodule ${}_A M_A$ is $G$-graded if $M = \oplus_{g \in G} M_g$ such that $A_h M_g \subseteq M_{hg}$ and $M_g A_h \subseteq M_{gh}$ for all $g,h \in G$.  We let $\Bimod{A}$ denote the category of $(A,A)$-bimodules and $\Bimod[G]{A}$ denote the category of $G$-graded bimodules and degree-preserving morphisms.  Both are abelian categories, but unlike usual categories of graded modules, $\Bimod[G]{A}$ does not appear to be equivalent to a module category when $G$ is nonabelian.  Nevertheless, when $A$ has a radical grading, graded bimodules admit graded projective covers, and these can be constructed by placing suitable gradings on the (ungraded) projective covers: if ${}_A M_A$ is graded, then its bimodule top $M/(JM + MJ)$ is graded and this defines a unique grading on the projective cover $P$ of $M$ such that the map $P \rightarrow M$ preserves degrees.  

In \cite{HCFDA}, Happel describes the terms in the minimal projective resolution $P^{\bullet}$ of $A$ as an $(A,A)$-bimodule.  The $r^{th}$ term is the projective bimodule $P^r = \oplus_{i,j} (Ae_i \otimes e_jA)^{m_{ij}}$ where $m_{ij} = \dim_k \Ext{A}{r}{S_i}{S_j}$.  This follows easily once one notices that $S_i \otimes_A P^{\bullet}$ is a minimal projective resolution of $S_i$.  Assuming that $J = J_G$, it follows from the above remarks that each $P^{r}$ can be graded to give a minimal graded projective resolution of the graded bimodule ${}_A A_A$.  For later reference, the grading on a summand of the form $Ae_i \otimes_k e_j A$ is obtained by letting $\deg(e_i \otimes e_j) = d$ for some $ d \in G$, which yields $$(Ae_i \otimes e_j A)_g = \bigoplus_{s,t \in G,\ sdt = g} A_s e_i \otimes e_j A_t$$ for each $g \in G$.  We denote this graded projective bimodule as $Ae_i \otimes e_jA[d]$.  Now considering $P^{\bullet}$ with this grading, we see that $S_i[d] \otimes_A P^{\bullet}_A$ gives a minimal graded projective resolution for $S_i[d]$.  Consequently, as a graded bimodule $P^r = \oplus_{i, j} \oplus_{d \in G} (Ae_i \otimes e_j A [d])^{m^d_{ij}}$ where $m^d_{ij} = \dim_k \Ext{A}{r}{S_i}{S_j[d]}$ in $\fmod[G]{A}$.  In particular, when $\Omega^r(S_i) \cong S_i$ for all $i$, we have $\Omega^r_{A^e}(A)$ generated in degree $e$ if and only if $\Omega^r(S_i) \cong S_i$ as graded modules for all $i$.

When the grading on $A$ is induced by assigning weights to the arrows of $Q$, Green shows that the category $\fmod[G]{A}$ is equivalent to $\fmod{B}$ where the quiver and relations for $B$ are a covering of those for $A$.  In this situation, $B$ is said to be a Galois cover of $A$.  In fact, such an equivalence exists for much more general gradings, and we employ the language of smash products to give an explicit description of the resulting algebra $B$ \cite{CoMo}.  For simplicity, we assume for the remainder of this article that $G$ is finite.  For $a \in A$ and $g \in G$, we write $a_g$ for the degree-$g$ component of $a$, and $p_g$ for the function $G \rightarrow k$ that sends $h$ to $\delta_{h,g}$.

\begin{defin} The {\bf smash product} of $A$ with $G$ is the $k$-algebra $A \# k[G]^* = \oplus_{g \in G} A p_g$ with multiplication given by $$ap_g \cdot bp_h = ab_{gh^{-1}}p_h, \ \ \forall a,b \in A.$$
\end{defin}

\noindent
{\it Remark.}  Even if $G$ is infinite, the smash product construction for categories can be applied to the $G$-graded category $\ind (\proj{A})$ of indecomposable projective $A$-modules to get a locally finite-dimensional category $\mathcal{B}$ as in \cite{CiMa}.  We also note that $\mathcal{B}$ is then a Galois cover of $\ind (\proj{A})$ with group $G$.  \\

We let $B = A \# k[G]^*$.  As each $e_i$ is homogeneous, a complete set of pairwise orthogonal primitive idempotents of $B$ is given by $\{e_i p_g\ |\ 1 \leq i \leq n,\ g \in G\}$.  We also observe that $\{p_h\ |\ h \in G\}$ is a set of pairwise orthogonal idempotents in $B$ that sum to $1$.  There is a free right $G$-action on this set, given by $p_h \cdot g = p_{hg}$, and this induces a right action of $G$ on $B$.  There is a natural embedding of algebras $i : A \rightarrow B$, sending $a \in A$ to $a \cdot 1 = \sum_{h \in G} a p_h \in B$, and one easily checks that $i$ identifies $A$ with the invariant subring $B^G$.  

As mentioned above, there is an isomorphism of categories $\Mod{B} \cong \Mod[G]{A}$ \cite{CoMo, CiMa, Green}, and we will often use it to identify graded $A$-modules with $B$-modules.  We also have a {\it pull-up functor} $F = - \otimes_A B_B : \Mod{A} \rightarrow \Mod{B}$, which is exact since ${}_A B$ is free.  If $M_A$ is a right $A$-module, then  $F(M) = M \otimes_A B = \oplus_{h \in G} M p_h$.  In terms of this decomposition, the right $B$-module structure on $F(M)$ is given by $(mp_h)(ap_g) = ma_{hg^{-1}} p_g$.  
The {\it push-down functor} $ i^* : \Mod{B} \rightarrow \Mod{A}$ is induced by the embedding $i : A \rightarrow B$.  Since this is essentially a restriction functor, it is exact and right adjoint to $F$.    As $G$ is finite, these functors restrict to the full subcategories of finitely-generated modules.

Before stating one last result on smash products, we review the definition of twisted bimodules.   If $\sigma \in \aut(A)$ is a $k$-algebra automorphism, and ${}_A M_A$ is an $(A,A)$-bimodule, we will write ${}_1M_{\sigma}$ for the twisted bimodule, where the left action of $A$ is the same as on $M$ but the right action of $A$ is twisted by $\sigma$: $m \cdot a = m a^{\sigma}$.  Equivalently, we have ${}_1 M_{\sigma} \cong {}_1 M \otimes_A {}_1 A_{\sigma}$.  Concerning twisted bimodules of the form ${}_1 A_{\sigma}$, we have the following simple observations: (1) $\sigma^{-1} : {}_1A_{\sigma} \stackrel{\cong}{\rightarrow} {}_{\sigma^{-1}}A_1$ is an isomorphism of bimodules;\ (2) ${}_1 A_{\sigma} \otimes_A {}_1 A_{\tau} \cong {}_1 A_{\tau \sigma}$; and (3) ${}_1A_{\sigma} \cong {}_1 A_1$ if and only if $\sigma$ is an inner automorphism.

\begin{lemma}  There is an isomorphism of $(B,B)$-bimodules $\displaystyle {}_B B \otimes_A B_B \cong \bigoplus_{x \in G} {}_1 B_x$.
\end{lemma}

\noindent
{\it Proof.} We have $${}_B B \otimes_A B_B  \cong  \bigoplus_{s \in G}  B \otimes p_s  \cong \bigoplus_{x \in G} \left( \bigoplus_{g \in G} Ap_g \otimes p_{gx^{-1}} \right).$$   Finally, the map sending $\sum_g a^g p_g \otimes p_{gx^{-1}} \mapsto \sum_g a^g p_g$ is easily seen to give a $(B,B)$-bimodule isomorphism $\oplus_{g \in G} Ap_g \otimes p_{gx^{-1}} \stackrel{\cong}{\rightarrow} {}_1 B_x$ for each $x \in G$.  $\Box$



\section{Lifting bimodules and resolutions}

Keeping the notation of the previous section, we now specify how to lift a graded $(A,A)$-bimodule to a $(B,B)$-bimodule, and apply this to construct a projective bimodule resolution for $B$ from one for $A$.  When $A$ is radical graded, we show that one of these resolutions is periodic if and only if the other is.  Still the definitions and basic properties proved below are valid for arbitrary group gradings, and we shall make use of them in greater generality in Section 6.


\begin{defin} Let ${}_A M_A$ be a $G$-graded bimodule, and fix $x \in G$.  We let $F_x(M)$ equal $F(M)$ as a right $B$-module (even as an $(A,B)$-bimodule), and we define a left $B$-module structure on $F(M)$ by the formula
$$ap_g \cdot m_kp_h = \left\{ \begin{array}{cl} am_kp_h, & \mathrm{if}\ g=khx \\ 0 & \mathrm{if}\ g\neq khx \end{array} . \right.$$ 
\end{defin}

We first check that this defines a left $B$-action on $F(M)$.  For $b p_l \in B$ and $a_i \in A_i$, we have \begin{eqnarray*} bp_l \cdot (a_ip_g \cdot m_kp_h) & = & \left\{ \begin{array}{cl} bp_l \cdot (a_i m_k p_h), & \mathrm{if}\ g=khx \\ 0 & \mathrm{if}\ g\neq khx \end{array} \right. \\ & = & \left\{ \begin{array}{cl} ba_i m_k p_h, & \mathrm{if}\ g=khx,\ l = ikhx \\ 0 & \mathrm{otherwise} \end{array} \right. ,\end{eqnarray*}
and 
\begin{eqnarray*} (bp_l \cdot a_i p_g) \cdot ( m_k p_h) & = & \left\{ \begin{array}{cl} (b a_i p_g) \cdot (m_kp_h), & \mathrm{if}\ i=lg^{-1} \\ 0 & \mathrm{if}\ i\neq lg^{-1} \end{array} \right. \\ & = & \left\{ \begin{array}{cl} b a_i m_k p_h, & \mathrm{if}\ l=ig,\ g=khx \\ 0 & \mathrm{otherwise} \end{array} \right. .
\end{eqnarray*}

We now check that the left and right $B$-actions on $F(M)$ commute.  We have 
\begin{eqnarray*} a p_g \cdot (m_k p_h \cdot b p_l) & = & a p_g \cdot (m_k b_{hl^{-1}} p_l) \\ & = & \left\{ \begin{array}{cl} a m_k b_{hl^{-1}} p_l, & \mathrm{if}\ g=khx \\ 0 & \mathrm{if}\ g\neq khx \end{array} \right. ,\end{eqnarray*}
and 
\begin{eqnarray*} (a p_g \cdot m_k p_h) \cdot bp_l & = & \left\{ \begin{array}{cl} (a m_k p_h) \cdot bp_l, & \mathrm{if}\ g = khx \\ 0 & \mathrm{if}\ g \neq khx \end{array} \right. \\ & = & \left\{ \begin{array}{cl} a m_k b_{hl^{-1}} p_l, & \mathrm{if}\ g=khx \\ 0 & \mathrm{if}\ g \neq khx \end{array} \right. .
\end{eqnarray*}

If $f : M \rightarrow N$ is a morphism in $\Bimod[G]{A}$ we can check that $F(f) : F(M) \rightarrow F(N)$ also respects the left $B$-action in this case, and so is a map of $(B,B)$-bimodules.  Hence $F_x : \Bimod[G]{A} \rightarrow \Bimod{B}$ is a functor.

We now establish several basic properties of these lifting functors.  Throughout, ${}_AM_A$ will denote a $G$-graded $(A,A)$-bimodule, and the bimodule ${}_A A_A$ is given the same grading as the algebra $A$.

\begin{lemma} $F_e({}_1 A_1) \cong {}_1 B_1$.
\end{lemma}

\noindent
{\it Proof.} Clearly, the natural isomorphism $F_e({}_1 A_1) = A \otimes_A B \longrightarrow B$ is a morphism of right $B$-modules.  We check that it is also a left $B$-module morphism.  For any $a,b \in A$ and any $h,g, l \in G$, we have $bp_g \cdot a_l p_h = ba_l p_h$ if $g = lh$ and it is $0$ otherwise.  Similarly, the multiplication in $B$ yields $bp_g a_l p_h = ba_l p_h$ if $l = gh^{-1}$, which holds if and only if $g=lh$, and is $0$ otherwise. $\Box$ \\ 

\begin{lemma} If $x \in G$, $F_x(M) \cong {}_1F_e( M)_x \cong {}_{x^{-1}}F_e( M)_1$.  In particular, $F_x({}_1 A_1) \cong {}_1 B_x$.
\end{lemma}

\noindent
{\it Proof.}  Consider the bijective map $f : mp_h \mapsto mp_{hx}$ on $F(M) = M \otimes_A B$.  We check that this is a bimodule morphism from $F_x(M)$ to ${}_1F_e( M)_x$.  We have $f(ap_g \cdot m_kp_h) = f(am_kp_h) = am_kp_{hx}$ if $g = khx$ and otherwise it is $0$.  We also have $ap_g \cdot f(m_kp_h) = ap_g \cdot m_kp_{hx} = am_k p_{hx}$ if $g = khx$ or $0$ otherwise.  On the right side, we have $f(mp_h ap_g) = f(ma_{hg^{-1}}p_g) = ma_{hg^{-1}}p_{gx}$ and $f(mp_h) \cdot ap_g = mp_{hx} ap_{gx} = ma_{hg^{-1}}p_{gx}$.  The second isomorphism follows from the remarks in Section 2.
$\Box$ \\

\begin{lemma} For $x \in G$, the functor $F_x$ is exact and takes projectives to projectives.
\end{lemma}

\noindent
{\it Proof.}  Since the usual pull-up functor $F : \Bimod[G]{A} \rightarrow \Mod{B}$ is exact and equals the composite of $F_x$ and the forgetful functor $\Bimod{B} \rightarrow \Mod{B}$, $F_x$ must be exact.  Since $F$ is additive, it suffices to show that $F_x(A \otimes_k A)$ is projective for any $G$-grading of the projective bimodule ${}_A A \otimes_k A_A$.  Let ${}_B Q_B  = F_x(A \otimes_k A)$, which is isomorphic to $A \otimes_k B$ as an $(A,B)$-bimodule.  Thus $B \otimes_A B \otimes_B Q \cong B \otimes_A (A \otimes_k B) \cong B \otimes_k B$  is a projective bimodule.  On the other hand, by Lemma 2.2, $\displaystyle B \otimes_A B \otimes_B Q \cong \bigoplus_{x \in G} {}_x Q_1$.  Since $Q$ is clearly a summand of the latter, it is projective.  $\Box$ \\

\noindent
{\it Remark.}  For the grading of the indecomposable projective $(A,A)$-bimodule ${}_A P_A = Ae_i \otimes_k e_j A$ with $e_i \otimes e_j$ in degree $e$, one can check that there is an isomorphism $$F_x(P) \cong \bigoplus_{s \in G} Be_ip_{sx} \otimes_k e_jp_s B$$ of $(B,B)$-bimodules.\\

\begin{lemma} Let ${}_AM_A$ and ${}_AN_A$ be $G$-graded bimodules.  If $M \cong N$ as ungraded bimodules and $F_e(M)$ is an indecomposable $(B,B)$-bimodule, then $F_e(M) \cong F_x(N)$ for some $x \in G$.
\end{lemma}

\noindent
{\it Proof.}  Let $f : M \rightarrow N$ be an $(A,A)$-bimodule isomorphism.  Tensoring with $B$, over $A$, on both sides yields and isomorphism of $(B,B)$-bimodules $1_B \otimes f \otimes 1_B : B \otimes_A F_e(M) \rightarrow B \otimes_A F_e(N)$.  Since $F_e(M)$ is a $(B,B)$-bimodule, we have isomorphisms $$B \otimes_A F_e(M) \cong B \otimes_A B \otimes_B F_e(M) \cong \bigoplus_{x \in G} {}_x F_e(M)_1$$
by Lemma 2.2, and similarly for $F_e(N)$.  The indecomposability of $F_e(M)$, along with the Krull-Schmidt theorem, now implies that $F_e(M) \cong {}_{x^{-1}}F_e(N)_1 \cong F_x(N)$ for some $x \in G$. $\Box$ \\ 


If $M$ is a graded bimodule and $d \in Z(G)$, we can shift the grading of $M$ by $d$ to get another graded bimodule $M[d]$ with $M[d]_g = M_{d^{-1}g}$ for all $g \in G$.  Clearly $M[d] \cong M$ as ungraded bimodules.  In this case it is easy to see that we have a $(B,B)$-bimodule isomorphism $F_e(M[d]) \cong F_d(M)$.

\begin{lemma}  Suppose $\sigma \in \aut(A)$ is a degree preserving automorphism such that $F_e({}_{\sigma} A_1) \cong {}_1 B_1$ as $(B,B)$-bimodules (we give ${}_{\sigma}A_1$ the same grading as $A$).  Then $\sigma^{|G|}$ is an inner automorphism.
\end{lemma}

\noindent
{\it Proof.}  Suppose that an isomorphism $f: B \rightarrow F_e({}_{\sigma}A_1)$ is given by sending $1_B$ to $\sum_{g \in G} m^g p_g$ for $m^g \in A$.  We first claim that each $m^g$ is a unit concentrated in degree $e$.  On one hand, we have $f(p_g) = f(1)p_g = m^g p_g$, and on the other $f(p_g) = p_g f(1) = p_g \sum_{h \in G} m^h p_h = \sum_{h \in G} m^h_{gh^{-1}} p_h$.  It follows that $m^g_h = 0$ if $h \neq e$ and $m^g_e = m^g$.    Now consider $f(ap_g) = \sum_{h \in G} a p_g m^hp_h = a^{\sigma} m^g p_g$.  Since $f$ is injective, $a^{\sigma} m^g \neq 0$ for all nonzero $a \in A$.  Surjectivity of $f$ now implies that $Am^g = A$, and hence $m^g$ must be a unit.  We also have $f(ap_g) = \sum_{h \in G} m^hp_h a p_g = \sum_{h \in G} m^h a_{hg^{-1}} p_g$.  Therefore, for all $h, g \in G$ and all $a_h \in A_h$ we have $$a_h^{\sigma}  = m^{hg} a_h (m^g)^{-1}.$$ 
We now apply this identity to show that the $m^{g}$ commute with each other.  Since $m^g \in A_e$ for all $g \in G$, we have $(m^g)^{\sigma} = m^g m^g (m^g)^{-1} = m^g$ and $(m^g)^{\sigma} = m^h m^g (m^h)^{-1}$ for any other $h \in G$.  We claim that $\sigma^{|G|}$ is conjugation by $m = \prod_{g \in G} m^g$.   Let $a_h \in A_h$ where $h$ has order $r$.  Then for any $x \in G$, $(a_h)^{\sigma^r} = m^{h^rx} \cdots m^{hx} a_h (m^x)^{-1} \cdots (m^{h^{r-1}x})^{-1}$, which is conjugation by $\prod_{g \in \gen{h}x} m^g$.  Thus, if we repeat as $x$ runs through a right transversal to $\gen{h}$ in $G$, we see that $(a_h)^{\sigma^{|G|}} = m a_h m^{-1}$.  $\Box$ \\

We now assume that $A$ has a radical grading and show that $A$ is periodic if and only if $B$ is.   For convenience, we assume that both $A$ and $B$ are indecomposable.  This ensures that the bimodules ${}_A A_A$ and ${}_B B_B$ are indecomposable.  (For gradings associated to Galois covers, we know that $B$ is indecomposable if and only if the grading on $A$ is {\it connected} in the terminology of \cite{Green}, i.e., if for all vertices $u, v  \in Q_0$ and each $g \in G$ there is an (undirected) walk in $Q$ from $u$ to $v$ of degree $g$ in $kQ$.)  As in the previous section, we let $P^{\bullet} : \cdots \rightarrow P_1 \rightarrow P_0 \rightarrow {}_A A_A \rightarrow 0$ be a minimal (graded) projective resolution of $A$ as an $(A,A)$-bimodule.  

\begin{therm} Suppose that $A$ is a $G$-graded $k$-algebra with homogeneous radical such that $B = A \# k[G]^*$ is indecomposable.  Then $A$ has a periodic projective resolution over $A^e$ if and only if $B$ has a periodic projective resolution over $B^e$.  Furthermore, if $p_A$ and $p_B$ denote the periods of $A$ and $B$ respectively, then $p_B \mid p_A \exp(G)$ and $p_A \mid p_B |G|$.
\end{therm}

\noindent
{\it Proof.}  First suppose that $A$ is periodic.  By lemmas 3.2 and 3.4, $F_e(P^{\bullet})$ will be a projective resolution of $B$.  Since $\Omega^r_{A^e} (A) \cong A$, Lemmas 3.5 and 3.3 show that $\Omega^r_{B^e}(B) \cong {}_1 B_x$ for some $x \in G$.  If $x^m = e$, then $\Omega^{rm}_{B^e}(B) \cong {}_1 B_{x^m} \cong B$.  Conversely, if $\Omega^r_{B^e}(B) \cong B$, then $\Omega^r_B(S) \cong S$ for every simple $B$-module $S$.  Equivalently, $\Omega^r_A (S_i) \cong S_i$ as graded modules for each simple $A$-module $S_i$.  As remarked in the previous section, this implies that $\Omega^r_{A^e}(A)$ is generated in degree $e$.  By Theorem 1.4 of \cite{GrSnSo}, $\Omega^r_{A^e}(A) \cong {}_1 A_{\sigma}$ for some automorphism $\sigma$ of $A$.  Moreover, in our case, the proof of this theorem easily yields that $\sigma$ preserves the grading on $A$ and this bimodule isomorphism is degree-preserving (where ${}_1 A_{\sigma}$ has the same grading as $A$).  Comparing the projective resolution $F_e(P^{\bullet})$ to a minimal projective resolution of $B$, we see that $F_e({}_1 A_{\sigma}) \cong \Omega^r_{B^e}(B) \cong B$.  Thus, by Lemma 3.6 $\sigma^{|G|}$ is inner, and hence $\Omega^{r|G|}_{A^e}(A) \cong {}_1 A_{\sigma^{|G|}} \cong A$.  $\Box$ \\

It would be interesting to determine whether we always have $p_A \mid p_B$.  Such a relation appears plausible and would make the computation of the periods of standard self-injective algebras of finite type significantly more tractable.  However, the following example shows that an automorphism $\sigma$ satisfying the hypotheses of Lemma 3.6 is not necessarily inner.  We let $A = P(\mathbb{L}_n)$ be the preprojective algebra associated to the generalized Dynkin graph $\mathbb{L}_n$ \cite{DPA}, i.e., it has quiver and relations
$$\xymatrix{0 \ar@(dl,ul)^{\epsilon = \bar{\epsilon}} \ar[r]<0.5ex>^{a_0} & 1  \ar[r]<0.5ex>^{a_1}  \ar[l]<0.5ex>^{\bar{a}_0} & 2  \ar[l]<0.5ex>^{\bar{a}_1} \ar@{}[r]|(0.4)\cdots \ar@{}[r] & n-2 \ar[r]<0.5ex>^{a_{n-2}}  & n-1 \ar[l]<0.5ex>^(.4){\bar{a}_{n-2}}}, \ \ \ \ \sum_{s(\alpha) = u} \bar{\alpha} \alpha = 0 \ \ (0 \leq u \leq n-1),$$ 
where we adopt the convention that $\bar{\bar{\alpha}} = \alpha$.  We give $A$ the $\mathbb{Z}/\gen{2}$-grading induced by the path length grading.  It is easy to see that $B = A \# k[\mathbb{Z}/\gen{2}]^*$ is then the preprojective algebra associated to the Dynkin graph $\mathbb{A}_{2n}$ (see Section 5).  If the characteristic of $k$ is not $2$, then both algebras are periodic of period $6$ \cite{DPA}.  Moreover, we have $\Omega^3_{A^e}(A) \cong {}_1 A_{\sigma}$ where $\sigma$ is the automorphism of $A$ induced by multiplying all arrows by $-1$.  In the notation of the proof of Lemma 3.6, if we let $m^e = 1$ and $m^x = -1$, then we have a $(B,B)$-bimodule isomorphism $F_e({}_1 A_{\sigma}) \cong B$ when ${}_1 A_{\sigma}$ is given the same grading as $A$.  However, since $B$ has period $6$, we can conclude that as graded bimodules $\Omega^3_{A^e}(A) \cong {}_1 A_{\sigma}[x]$, and hence $\Omega^3_{B^e}(B) \cong {}_1 B_x$.  It is well-known that $x$ is in fact a Nakayama automorphism of $B=P(\mathbb{A}_{2n})$, induced by the reflection of $\mathbb{A}_{2n}$.  \\


Given the algebra $B = A \# k[G]^*$, we can recover $A$, up to Morita equivalence, as the skew group algebra $B*G$ \cite{CiMa, CoMo}.  Recall that $B*G$ is a free $B$-module on $G$ with multiplication given by  $ag \cdot b h = a b^{g^{-1}} gh$ where the action of $G$ on $B$ is as described following Definition 2.1.
Thus, rephrasing the above theorem for skew group algebras yields the following.

 \begin{coro}  Suppose the finite group $G$ acts via automorphisms on a basic $k$-algebra $B$ with a free action on a complete set of pairwise orthogonal primitive idempotents for $B$.  Then $B$ is periodic if and only if $B*G$ is periodic.  \end{coro}

\noindent
{\it Proof.} In \cite{CiMa}, Cibils and Marcos show that such a $B$ is isomorphic to the smash product of $\und{B*G}$ with $G$, where $\und{B*G}$ is a basic version of $B*G$.  We claim that the assumption that $B$ is basic forces the $G$-grading on $A = \und{B*G}$ to be a radical grading.  If it is not a radical grading, then we must have a strict inequality $J_G(A) \subset J(A)$ and $J(B) = \sum_{g \in G} J_G(A)p_g$  by \cite{CoMo}.  Thus $$\dim_k B/J(B) > \dim_k B/J(A)B = \dim_k (A/J(A) \otimes_A B) = |G| \cdot \dim_k A/J(A),$$ which equals the cardinality of a complete set of pairwise orthogonal primitive idempotents of $B = A \#k[G]^*$.  But this would contradict the assumption that $B$ is basic.  $\Box$ \\


In \cite{CiRe}, Cibils and Redondo establish a spectral sequence for Hochschild cohomology associated to a Galois covering $B \rightarrow A$.  Our approach provides information on how the minimal projective resolutions for $A$ and $B$ are related in this case, and could perhaps be applied to give a more direct comparison between the Hochschild (co)homology of the two algebras.

\section{Stable Auslander Algebras}

We now assume that $A$ is a self-injective algebra of finite representation type in order to compare periodicity properties of $A$ and its stable Auslander algebra.  This problem is explored in greater generality in \cite{Buch}, and we now review Buchweitz's results in our simplified context.  We let $M_{A}$ be a (basic) representation generator for $A$, i.e., $M$ is the direct sum of one representative from each isomorphism class of indecomposable right $A$-modules, and let $M'$ denote the direct sum of all nonprojective indecomposable summands of $M$.  We let $\Lambda = \End{A}{M}$ be the Auslander algebra of $A$, and we let $\Gamma = \stEnd{A}{M} \cong \stEnd{A}{M'}$ be the stable Auslander algebra of $A$.  Note that $\Gamma$ is just the quotient of $\Lambda$ by the ideal $\Lambda \pi \Lambda$ where $\pi \in \End{A}{M}$ denotes the projection from $M$ onto $A$.  It follows from results of \cite{SEAA}, and is proved directly in \cite{Buch}, that $\Gamma$ is also self-injective.  The functor $\Hom{A}{M}{-} : \fmod{A} \rightarrow \proj{\Lambda}$ is an equivalence and induces an equivalence $\stHom{A}{\und{M}}{-} : \stmod{A} \rightarrow \proj{\Gamma}$.  It follows that the quiver of $\Lambda$ is the AR-quiver of $A$, and the quiver of $\Gamma$ is the stable AR-quiver of $A$.  Furthermore, if $A$ is standard, then the relations for $\Lambda$ and $\Gamma$ are just the mesh relations associated with these translation quivers.

As shown in sections 5 and 6 of \cite{Buch}, the natural ring homomorphism $\Lambda \rightarrow \Gamma$ is pseudoflat and $L = \Tor{\Lambda}{2}{\Gamma}{\Gamma} \cong \Omega^3_{\Gamma^e}(\Gamma)$ as $(\Gamma, \Gamma)$-bimodules.  Thus, tensoring with $L$ induces $\Omega^3_{\Gamma}$ on $\stmod{\Gamma}$.  Furthermore, we have $(\Gamma, \Gamma)$-bimodule isomorphisms $L^{\otimes i} \cong \stHom{A}{M}{\Omega^i M} \cong \stHom{A}{M'}{\Omega^i M'}$ for $i \geq 0$.  The left $\Gamma$-module structure on $\Omega^i M'$ is given via an isomorphism $\Omega^i M' \cong M'$, which exists since $M$ is a representation generator and $A$ is self-injective.

\begin{propos} The following are equivalent for an integer $n \geq 1$: 
\begin{enumerate}
\item $\Omega^{3n}_{\Gamma^e} (\Gamma) \cong  \Gamma$.
\item $\Omega^n M' \cong M'$ as $(\Gamma, A)$-bimodules.
\item There is an isomorphism $\Omega^n_{A} \cong 1_{\stmod{A}}$ of functors on $\stmod{A}$.
\end{enumerate}
\end{propos}

\noindent
{\it Proof.}  The equivalence of (2) and (3) follows from the definition $\Gamma = \stEnd{A}{M'}$, while $(2) \Rightarrow (1)$ follows from the isomorphisms cited above.  For $(1) \Rightarrow (2)$, the isomorphism $\Omega^{3n}_{\Gamma^e} (\Gamma) \cong L^{\otimes n} \cong  \Gamma$ yields an isomorphism $\xi : \stHom{A}{M'}{\Omega^n M'} \rightarrow \stHom{A}{M'}{M'}$ of $(\Gamma, \Gamma)$-bimodules.  Now let $\varphi = \xi^{-1}(1_{M'}) : M' \rightarrow \Omega^n M'$.  Then $\varphi$ is a $(\Gamma, A)$-bimodule homomorphism.  By Yoneda's lemma $\xi$ is induced by a morphism $\chi: \Omega^n M' \rightarrow M'$ in $\stmod{A}$, and since $\chi \varphi = 1_{M'}$, $\varphi$ must be an isomorphism.  $\Box$ \\

In order to state the main result of this section we need a simple definition.  We say that $A$ is {\it Schurian} if $\dim_k e_i A e_j \leq 1$ for all $e_i, e_j$ belonging to a complete set of pairwise orthogonal primitive idempotents for $A$ (i.e., the entries of the Cartan matrix of $A$ are $0$ or $1$). 

\begin{therm}  Let $A$ be a basic, indecomposable self-injective $k$-algebra of finite representation type, and $\Gamma$ its stable Auslander algebra.  
\begin{enumerate}
\item If $\Omega^n_{A^e}(A) \cong A$, then $\Omega^{3n}_{\Gamma^e}(\Gamma) \cong \Gamma$.
\item If $A$ is Schurian and $\Omega^{3n}_{\Gamma^e}(\Gamma) \cong \Gamma$, then $\Omega^{n}_{A^e}(A) \cong A$.
\end{enumerate}
\end{therm}

\noindent
{\it Remark.}  We note that the period of $\Gamma$ is divisible by $3$ as long as the tree class of $A$ is not $\mathbb{A}_1$ or $\mathbb{A}_2$.  In this case, $A$ has an almost split sequences with at least 2 indecomposable nonprojective summands in the middle term, and hence the projective resolution for the corresponding simple $\Gamma$-module (cf. section I.3 of \cite{SEAA}) has a decomposable $n^{th}$ term if and only if $n \equiv 1 \mod 3$.  On the other hand, since the preprojective algebra $P(\mathbb{A}_2)$ has period $2$ \cite{ErdSna2}, the results of the previous section imply that the period of the mesh algebra of $\mathbb{Z}\mathbb{A}_2/\gen{\tau^m}$, which is an $m$-fold covering of $P(\mathbb{A}_2)$ (see below), is not divisible by $3$ whenever $3 \nmid m$. \\

The first part of the theorem follows from the proposition and the comments preceding it.  For the second part, we need to investigate automorphisms of Schurian algebras and the corresponding twisted bimodules.  We thus suppose that $A$ is Schurian and that $\sigma \in \aut(A)$ is an automorphism fixing $e_i$ for each $i$.  We fix a presentation of $A$ as a path algebra of a quiver $Q = (Q_0, Q_1)$ with relations, where the vertex set $Q_0$ of the quiver is identified with the complete set of primitive idempotents $\{e_i\}_{1\leq i \leq n}$ and the arrows correspond to chosen elements of $e_j J_{A} e_i / e_j J_{A}^2 e_i$.  For each arrow $\alpha$ of $Q$, there exists $c_{\alpha} \in k^*$ such that $\sigma(\alpha) = c_{\alpha} \alpha$.

\begin{lemma} Let $\sigma \in \aut(A)$ be as above.  Then ${}_1 A_{\sigma} \cong {}_1 A_1$ if and only if there exist $d_i \in k^*$ for each $i \in Q_0$ such that $c_{\alpha} = d_j/d_i$ for all $i,j \in Q_0$ and all arrows $\alpha$ from $i$ to $j$. 
\end{lemma}

\noindent
{\it Proof.}  For the forward direction, assume that $f : {}_1 A_{\sigma} \rightarrow {}_1 A_1$ is a bimodule isomorphism.  We have $f(e_i) = e_if(e_i)e_i$, which implies that $f(e_i) = d_i e_i$ for some $d_i \in k^*$.  If $\alpha$ is an arrow from $i$ to $j$, we have $\alpha = e_j \alpha e_i$, and thus $f(\alpha) = f(e_j \alpha) = f(e_j) c_{\alpha}^{-1} \alpha = c_{\alpha}^{-1} d_j \alpha$.  Similarly, we get $f(\alpha) = f(\alpha e_i) = \alpha f(e_i) = d_i \alpha$.  This shows that $c_{\alpha} = d_j/d_i$.

Conversely, suppose that nonzero scalars $d_i$ exist so that $c_{\alpha} = d_j/d_i$ for all $i, j \in Q_0$ and all arrows $\alpha$ from $i$ to $j$.  Then $\alpha^{\sigma} = d_j \alpha d_i^{-1} = u \alpha u^{-1}$, where $u = \sum_{i \in Q_0} d_i e_i$ is a unit.  Hence $\sigma$ is an inner automorphism and the result follows. $\Box$ \\

\begin{lemma} Suppose $A$ is self-injective and Schurian, and $\sigma \in \aut(A)$ fixes $e_i$ for all $i$.  If $-\otimes_{A} A_{\sigma} : \stmod{A} \rightarrow \stmod{A}$ is isomorphic to $Id_{\stmod{A}}$, then ${}_1 A_{\sigma} \cong {}_1 A_1$ as bimodules.
\end{lemma}

\noindent
{\it Proof.}  Let $\eta : -\otimes_{A} A_{\sigma} \rightarrow Id_{\stmod{A}}$ be an isomorphism.  For any indecomposable, nonprojective $A$-module $M_{A}$, $\eta_M : M_{\sigma} = M \otimes_{A} {}_1 A_{\sigma} \rightarrow M$ is an isomorphism in the stable category.  Thus any lift of $\eta_M$ to $\fmod{A}$ is an isomorphism, and we fix such a lift for each $M$ and continue to denote these as $\eta_M$.  

Now consider two arrows $\alpha$ and $\beta$ from $j$ to $l$ and from $i$ to $j$, respectively, such that $\alpha \beta \neq 0$.  We consider the surjective map $f_{\alpha} : \beta A \rightarrow \alpha \beta A$ between indecomposable nonprojective $A$-modules, which is given by left-multiplication by $\alpha$, and the inclusion $g_{\beta} : \alpha \beta A \rightarrow \alpha A$.  Abbreviating $\eta_{\alpha A}$ as $\eta_{\alpha}$ etc., we have diagrams which are commutative in the stable category:
$$\xymatrix{\beta A_{\sigma} \ar[r]^{f_{\alpha}} \ar[d]_{\eta_{\beta}}^{\cong} & \alpha \beta A_{\sigma} \ar[r]^{g_{\beta}} \ar[d]_{\eta_{\alpha \beta}}^{\cong} & \alpha A_{\sigma} \ar[d]_{\eta_{\alpha}}^{\cong} \\ \beta A \ar[r]_{f_{\alpha}} & \alpha \beta A \ar[r]_{g_{\beta}} & \alpha A}$$
Since $A$ is Schurian and $\beta A$ has a simple top, $\eta_{\beta}$ must be multiplication by some nonzero scalar $c^{\beta}$.  Defining $c^{\alpha \beta}$ and $c^{\alpha}$ similarly, the commutativity of the left square shows that $c^{\alpha \beta} = c^{\beta}$, since otherwise the difference of the two maps would be surjective and thus could not factor through a projective.  The same reasoning shows that $c^{\beta} = c^{p}$ for any nonzero path $p$ starting at $i$.  Similarly, the commutativity of the right square shows $c^{\beta} = c^{\alpha}c_{\beta}^{-1}$, otherwise the difference of the two maps would be injective and hence could not factor through a projective.  We now have $c_{\beta} = c^{\alpha}/c^{\beta} = d_j/d_i$ where we let $d_i = c^{\beta}$, which only depends on $i$ as noted above, and $d_j = c^{\alpha}$, which only depends on $j$.  From the previous lemma, we conclude that ${}_1 A_{\sigma} \cong {}_1 A_1$ as bimodules. $\Box$\\

\noindent
{\it Remark.}  We do not know whether the assumptions on $A$ in the above lemma are truly necessary. \\


\noindent
{\it Proof of Theorem 4.2.}  It remains to prove (2).  According to Proposition 4.1, we have an isomorphism of functors $\Omega^n_{A} \cong Id_{\stmod{A}}$ on $\stmod{A}$.  Since we have $\Omega^n(S) \cong S$ for all simples $S$, Theorem 1.4 in \cite{GrSnSo} shows that $\Omega^n_{A^e} (A) \cong {}_1 A_{\sigma}$ for some $\sigma \in \aut(A)$, which fixes each $e_i$.  Finally, since $- \otimes_{A} {}_1 A_{\sigma} \cong - \otimes_{A} \Omega^n_{A^e}(A) \cong \Omega^n_{A} \cong Id_{\stmod{A}}$, Lemma 4.4 implies that $\Omega^n_{A^e}(A) \cong A$.  $\Box$ \\

 
 
\section{Self-injective algebras of finite representation type}

We now combine the results of the previous two sections to show that any standard self-injective algebra $A$ of finite representation type has a periodic projective resolution over its enveloping algebra $A^e$.  Not only does this imply that the Hochschild (co)homology groups of such an algebra are periodic, but by the results of \cite{GrSnSo} we can conclude that the Hochschild cohomology ring modulo the ideal generated by homogeneous nilpotent elements is isomorphic to $k[x]$ where the degree of $x$ is the period of the resolution.  As usual, all algebras we consider are assumed to be basic and indecomposable.

We briefly review the definitions of preprojective algebras first.  If $\Delta$ is a Dynkin graph, we can obtain a quiver $Q_{\Delta}$ by replacing each edge with a pair of arrows $\alpha$ and $\overline{\alpha}$ in opposite directions.  We regard $\alpha \mapsto \overline{\alpha}$ as an involution on the arrows of $Q_{\Delta}$ with $\overline{\overline{\alpha}} = \alpha$.  The preprojective algebra $P(\Delta)$ is then defined to be the path algebra modulo relations $k Q_{\Delta}/I$ where $I$ is the ideal generated by the sums $\sum_{s(\alpha) = u} \overline{\alpha} \alpha$ which range over all arrows $\alpha$ of $Q_{\Delta}$ with source $s(\alpha) = u$, for each vertex $u$ of $Q_{\Delta}$.  For any integer $m>1$, we can give $P(\Delta)$ a $\mathbb{Z}/\gen{m}$-grading by assigning to each pair of arrows $\alpha, \overline{\alpha}$ degrees $0$ and $1$.  The smash product $P(\Delta) \#k[\mathbb{Z}/\gen{m}]^*$ is then easily seen to be isomorphic to the mesh algebra associated to the translation quiver $\mathbb{Z}\Delta/\gen{\tau^m}$.  Since $P(\Delta)$ is known to be periodic, so are these mesh algebras by Theorem 3.7.  In fact, they were originally shown to be periodic in \cite{PAAK} using the fact that they are almost Koszul.

\begin{therm}[cf. 3.10 in \cite{ErdSko2}] Any standard self-injective algebra $A$ of finite representation type is periodic.
\end{therm}

\noindent
{\it Proof.}  In \cite{UCQR}, Mart\'{i}nez-Villa and de la Pe\~{n}a prove that any basic, indecomposable standard algebra $A$ of finite representation type admits a finite Galois cover $p : B \rightarrow A$ with $B$ Schurian.  Hence, by Theorem 3.7, we may assume that $A$ is Schurian.  If $\Gamma$ is the stable Auslander algebra of $A$, by Theorem 4.2 it suffices to prove that $\Gamma$ is periodic.  Since $A$ is standard, $\Gamma$ is isomorphic to the path algebra of the AR-quiver of $A$ modulo the ideal of mesh relations $I$.  By Riedtmann's structure theorem \cite{ADK}, the AR-quiver of $A$ has the form $\mathbb{Z}\Delta/\gen{\zeta \tau^{-r}}$, where $\Delta$ is the tree class of $A$ (an oriented Dynkin graph), $\zeta$ is an admissible automorphism of $\mathbb{Z}\Delta$ of finite order $t$, and $\tau$ is the translation.  Moreover, as $\zeta$ commutes with $\tau$ and $\gen{\zeta \tau^{-r}}$ is infinite cyclic, it follows that there is a finite Galois covering $\tilde{\Gamma} = k(\mathbb{Z}\Delta/\gen{\tau^{rt}})/I' \rightarrow \Gamma$, where $I'$ is the ideal of mesh relations for the new quiver.  As noted above, we also have a Galois $\mathbb{Z}/\gen{rt}$-covering $\tilde{\Gamma} \rightarrow P(\Delta)$ of the preprojective algebra of $\Delta$, which is obtained by factoring out the group of automorphisms of $\tilde{\Gamma}$ generated by $\tau$.  Since $P(\Delta)$ is periodic with period dividing $6$ (this was originally proved by Schofield, but see \cite{ErdSna2}), Theorem 3.7 implies that $\tilde{\Gamma}$ and $\Gamma$ are as well.  $\Box$ \\

We conclude this section by applying the strategy of the above proof to calculate upper bounds for the periods of the standard self-injective algebras of finite type.  Since the period of such an algebra is invariant under derived equivalence (cf. 2.2 in \cite{PAAK}), it suffices to look at one representative algebra from each derived equivalence class.  Such a list is given in \cite{DECSA}, and we refer the reader to the appendix of \cite{Asa2} for presentations of these algebras by quivers and relations.  These algebras are distinguished by their {\it type} $(\Delta, f, t)$, where $\Delta$ is the tree-class, $f = r / m_{\Delta}$ is the frequency, and $t$ is the order of $\zeta$ as in the above proof.  Here $m_{\Delta}$ equals $n, (2n-3), 11, 17$ or $29$ when $\Delta$ is $\mathbb{A}_n, \mathbb{D}_n, \mathbb{E}_6, \mathbb{E}_7$ or $\mathbb{E}_8$ respectively.  Note that the Coxeter number of $\Delta$ is $h_{\Delta} = m_{\Delta} + 1$.  We compile our results in the table at the end of the section.  The groups occuring below are all cyclic and we shall break with our previous notation to write them additively.

The precise periods are already known for algebras of tree class $\mathbb{A}_n$ \cite{EH, EHS} (see Table 5.2).  In order to bound the periods of the remaining algebras with tree class $\mathbb{D}$ or $\mathbb{E}$, we first focus on the mesh algebras $\Gamma$ of translation quivers of the form $\mathbb{Z}\Delta/\gen{\tau^m}$.  As noted above, such a $\Gamma$ is isomorphic to the smash product of the preprojective algebra $P(\Delta)$ with $\mathbb{Z}/\gen{m}$, where the grading on $P(\Delta)$ is given by assigning degrees $0$ and $1$ to each pair of arrows associated to an edge of $\Delta$.  With respect to the usual path-length grading, $P(\Delta)$ is $(h_{\Delta}-2, 2)$-Koszul, and thus $\Omega^3(P(\Delta))$ is generated in degree $h_{\Delta}$ and $\Omega^6(P(\Delta)) \cong P(\Delta)[2h_{\Delta}]$ as graded bimodules \cite{PAAK}.  With respect to our ``half-grading'', however, it is not difficult to see that $\Omega^6(P(\Delta)) \cong P(\Delta)[h_{\Delta}]$. Thus, by the remarks following Lemma 3.5, $\Omega^6(\Gamma) \cong {}_1 \Gamma_{h_{\Delta}}$, and it follows that $p_{\Gamma} | 6m/(h_{\Delta}, m)$.  To see that equality holds (assuming $\chr(k) \neq 2$), we look at two cases.  First, if $\Delta$ is $\mathbb{D}_n$ with $n$ odd or $\mathbb{E}_6$, then the sixth syzygy of $P(\Delta)$ is the first to fix all simple $P(\Delta)$-modules.  It follows that $6m/(h_{\Delta}, m)$ is the smallest simultaneous period of all graded simple $P(\Delta)$-modules, and thus of all simple $\Gamma$-modules.  Hence $p_{\Gamma} = 6m/(h_{\Delta}, m)$.  Now suppose that $\Delta = \mathbb{D}_{2n}, \mathbb{E}_7$ or $\mathbb{E}_8$.  Here, we have $\Omega^3(S) \cong S[h_{\Delta}/2]$ for each simple $P(\Delta)$-module $S$, and thus $3m/(h_{\Delta}/2, m) | p_{\Gamma}$.  Notice that if $m$ is even, this agrees with our upper bound for $p_{\Gamma}$.  If $m$ is odd we have $p_{P(\Delta)} = 6 | m p_{\Gamma}$ by Theorem 3.7, and hence $2 | p_{\Gamma}$.  Thus, we again have $p_{\Gamma} = 6m/(h_{\Delta}, m)$.  

In case the characteristic is $2$ and $\Delta = \mathbb{D}_{2n}, \mathbb{E}_7$ or $\mathbb{E}_8$, $P(\Delta)$ has period $3$ \cite{ErdSna2}, and we have $\Omega^3(P(\Delta)) \cong P(\Delta)[h_{\Delta}/2]$.  Similar to before, the upper bound we obtain for $p_{\Gamma}$ is now $3m/(h_{\Delta}/2, m)$, which coincides with the lower bound computed above.  Thus $p_{\Gamma} = 3m/(h_{\Delta}/2, m)$ in these cases (note that this value differs from the previous one only when $m$ is odd).
 
If $A$ has type $(\Delta, f, 1)$, its stable Auslander algebra $\Gamma$ is the mesh algebra of $\mathbb{Z}\Delta/\gen{\tau^{m_{\Delta}f}}$.  Thus, provided $A$ is Schurian, its period will be $p_{\Gamma}/3$ by Theorem 4.2 (see Table 5.2 for precise values).  It is straightforward to check that the representative algebras are Schurian precisely when $f>1$.  If $f \leq 1$, then we have only found the period of the syzygy functor on $\stmod{A}$.  Upper bounds for the actual periods of non-Schurian algebras can of course be obtained by passing to a Schurian cover, in which case we have to multiply the upper bound for the Schurian algebra by the degree of the covering.  This appears to yield poor upper bounds in general and so we omit them.

\setcounter{table}{1}
\renewcommand{\thetable}{\thesection .\arabic{table}}
\renewcommand{\arraystretch}{1.5}

\begin{table}
\begin{tabular}{|c|c|c|} \hline
Type of $A$ & Additional Cases & Period $p_{A}$ \\ \hline \hline
$(\mathbb{A}_n, s/n, 1)$ & char $k = 2, n=1$ and $2 \nmid s$ &   $ s$ \\ \cline{2-3}
			& otherwise & $ \frac{2s}{(s, n+1)}$ \\ \hline \hline
$(\mathbb{A}_{2m+1}, s, 2)$ &  char $k = 2$ and $2| \frac{s+m+1}{(s, m+1)}$ & $  \frac{s(2m+1)}{(s, m+1)}$ \\ \cline{2-3} 
			       & otherwise & $  \frac{2s(2m+1)}{(s, m+1)}$ \\ \hline \hline
$(\mathbb{D}_n, s, 1)$ & char $k=2$, $2 | n$ and $2 \nmid s$ & $\frac{s(2n-3)}{(s, 2n-2)}$ \\ \cline{2-3} 
			       & otherwise & $ \frac{2s(2n-3)}{(s, 2n-2)}$ \\ \hline \hline
$(\mathbb{D}_n, s, 2)$ & $2 \nmid n$ and $2 \mid \frac{s + n-1}{(s, n-1)}$ & $ \frac{s(2n-3)}{(s, n-1)}, \ \  \frac{2s(2n-3)}{(s, n-1)}$ or  $\frac{4s(2n-3)}{(s, n-1)}$ \\ \cline{2-3} 
			       & otherwise & $\ \frac{2s(2n-3)}{(s, n-1)}$ or $\frac{4s(2n-3)}{(s, n-1)}$ \\ \hline \hline
$(\mathbb{D}_{3m}, s/3, 1)$ & char $k=2$, $2|m$ and $2 \nmid s$ & $ \frac{s(2m-1)}{(s, 6m-2)}$ \\ \cline{2-3} 
			       & otherwise & $ \frac{2s(2m-1)}{(s, 6m-2)}$ \\ 
\hline \hline
$(\mathbb{D}_{3m}, 1/3, 1)$ & $A$ nonstandard & $\scriptstyle{2m-1,\ 2(2m-1)}$ or $\scriptstyle{4(2m-1)}$ \\
\hline \hline
$(\mathbb{D}_4, s, 3)$ & char $k=2$ and $2 \nmid s$ & $ 5s$ or  $15s$ \\ \cline{2-3} 
			       & otherwise & $ \frac{10s}{(s,2)}$ or $\frac{30s}{(s, 2)}$ \\ \hline \hline
$(\mathbb{E}_6, s, 1)$ & - & $  \frac{22s}{(s,12)}$ \\  \hline \hline
$(\mathbb{E}_6, s, 2)$ &  $ s \equiv 2 \mod 4$ & $ \frac{11s}{(s, 6)}, \ \ \frac{22s}{(s,6)}$ or $\frac{44s}{(s, 6)} $ \\  \cline{2-3} & $s \not \equiv 2 \mod 4$ & $ \frac{22s}{(s,6)}$ or $\frac{44s}{(s, 6)} $ \\  \hline \hline
$(\mathbb{E}_7, s, 1)$ & char $k=2$ and $2 \nmid s$ & $ \frac{17s}{(s,18)}$ \\ \cline{2-3} 
			       & otherwise & $  \frac{34s}{(s, 18)} $ \\ \hline \hline
$(\mathbb{E}_8, s, 1)$ & char $k=2$ and $2 \nmid s$ & $ \frac{29s}{(s,30)}$ \\ \cline{2-3} 
			       & otherwise & $\frac{58s}{(s,30)}$ \\ \hline
\end{tabular}
\caption{Periods of self-injective algebras of finite type in terms of their type $(\Delta, f, t)$.  For $f \leq 1$ and $\Delta \neq \mathbb{A}_n$ the values shown for the standard algebras only apply to the functorial period of $\Omega$.  When $t > 1$ or the algebra is nonstandard, we do not know if all possibilities occur.}  
\end{table}

We now consider those algebras of tree class $\mathbb{D}$ or $\mathbb{E}$ with torsion order $t >1$, which occur when the type of $A$ is $(\mathbb{D}_n, s, 2), (\mathbb{D}_4, s, 3)$ or $(\mathbb{E}_6, s, 2)$.  As in the proof of Theorem 5.1, the mesh algebra $\tilde{\Gamma}$ of the translation quiver $\mathbb{Z}\Delta/\gen{\tau^{m_{\Delta}ft}}$ is a $t$-fold cover of the stable Auslander algebra of $A$.  The period of $A$ is thus at most $t/3$ times the period of $\tilde{\Gamma}$, which was calculated above (see Table 5.2 for precise upper bounds).   Notice that for type $(\mathbb{D}_n, s, 2)$ this upper bound is independent of the characteristic since $m_{\Delta}ft$ is even.

In order to compute lower bounds, we treat each algebra separately and use functorial isomorphisms on their universal covers to deduce the order of $\Omega$ as a permutation on isomorphism classes of modules.  For type $(\mathbb{D}_n, s, 2)$, $\tau$ has order $2s(2n-3)$ and $\tau^{(2n-3)s}$ induces a permutation $\sigma$ of order $2$ on the indecomposable nonprojective $A$-modules (cf. proof of Prop. 2.5 in \cite{BiaSko}), while $\Omega$ coincides with $\tau^{n-1}$ if $n$ is even and with $\sigma \tau^{n-1}$ if $n$ is odd (cf. Prop. 4.2 in \cite{ErdSko}).  If $n$ is even, we see that $\Omega$ has order $2s(2n-3)/(s,n-1)$.  If $n$ is odd, $\Omega^{(2n-3)s/(s,n-1)} = \sigma^{((2n-3)s + n-1)/(s, n-1)}$, and thus $\Omega$ has order $s(2n-3)/(s,n-1)$ if $(s+n-1)/(s,n-1)$ is even, and it has order $2s(2n-3)/(s,n-1)$ otherwise.  

For type $(\mathbb{D}_4, s, 3)$, $\tau$ has order $15s$ and $\Omega = \tau^3$ on objects by Proposition 4.2 in \cite{ErdSko}.  Thus $\Omega$ has order $5s$.  However, when $\chr (k) \neq 2$, the functorial period of $\Omega$ must be even since the period of $P(\mathbb{D}_4)$ is even.  Hence, in this case, we get the lower bound $10s/(s,2)$.

Type $(\mathbb{E}_6, s, 2)$ is similar to $(\mathbb{D}_n, s, 2)$.  Here, $\tau$ has order $22s$, $\tau^{11s}$ induces a permutation $\sigma$ of order $2$, and $\Omega$ coincides with $\sigma \tau^6$ on objects (cf. Prop. 2.5 in \cite{BiaSko} and Prop. 4.2 in \cite{ErdSko}).  Thus $\Omega^{11s/(s,6)} = \sigma^{(11s+6)/(s,6)}$, and $\Omega$ has order $11s/(s,6)$ if $(s+6)/(s,6)$ is even (if and only if $s \equiv 2 \mod 4$) or order $22s/(s,6)$ otherwise. \\



Finally, recall that a self-injective algebra $A$ is said to be {\it stably $d$-Calabi-Yau} if there is an isomorphism of triangulated functors $\nu \cong \Omega^{-(d+1)}$ on $\stmod{A}$, where $\nu = - \otimes_{A} DA$ is the Nakayama equivalence.  If $A$ is symmetric, then $\nu \cong Id_{\stmod{A}}$ and hence $A$ is $d$-Calabi-Yau if and only if $d+1$ equals the order of $\Omega$ as a functor on $\stmod{A}$.  In particular, the algebra $A = P(\mathbb{L}_2)$ is a finite-type symmetric algebra of type $(\mathbb{D}_6, 1/3, 1)$.  If $\chr (k) \neq 2$ we find that the syzygy functor has order $6$, even though $\Omega^3(M) \cong M$ for every nonprojective indecomposable $A$-module $M$ \cite{DPA}.  It follows that $A$ is stably $5$-Calabi-Yau and not stably $2$-Calabi-Yau as claimed in \cite{ErdSko}.  In fact, the same is true for the algebras $P(\mathbb{L}_n)$, as it can be directly verified that $\Omega^3 \not \cong Id_{\stmod{A}}$ for these algebras using the description of $\Omega^3$ in \cite{DPA}.  Theorem 4.3 of \cite{ErdSko} includes similar errors for the symmetric algebras of tree classes $\mathbb{D}_{2n}, \mathbb{E}_7$ and $\mathbb{E}_8$ when $\chr (k)  \neq 2$.  As can be gleaned from Table 5.2, the stable Calabi-Yau dimension of the standard symmetric algebra with type $(\Delta, 1/r, 1)$ is $2m_{\Delta}/r -1$, or else $m_{\Delta}/r -1$ for $\Delta = \mathbb{A}_1, \mathbb{D}_{2n}, \mathbb{E}_7$ or $\mathbb{E}_8$ in characteristic $2$.  The error appears to arise from the (mistaken) assumption that an isomorphism of functors on the stable category of the universal cover of $A$ induces an isomorphism between the induced functors on the stable category of $A$.  For instance, for $A = P(\mathbb{L}_2)$ one has $\Omega \cong \tau^5$ over the universal cover of $A$, but not over $A$.  Earlier we saw that $P(\mathbb{A}_4)$ is a double cover of $A$, and it is stably $2$-Calabi-Yau since $\Omega^3$ is isomorphic to the Nakayama functor on $\stmod{P(\mathbb{A}_4)}$.  Thus we see that stable Calabi-Yau dimensions may indeed increase upon passage to the orbit algebra in a Galois cover.

\section{Nonstandard algebras}

Finally, we turn to the class of nonstandard indecomposable self-injective algebras of finite representation type.  These algebras arise only in characteristic $2$ as socle deformations of standard self-injective algebras of type $(\mathbb{D}_{3m}, 1/3, 1)$ for $m \geq 2$.  In particular, each has type $(\mathbb{D}_{3m}, 1/3,1)$ for some $m \geq 2$, and Asashiba has shown that two such are derived equivalent if and only if they have the same type \cite{DECSA}.  We thus focus on one representative algebra of each type, and these are given by the quivers 
$$
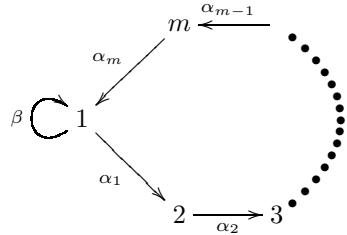
\xymatrix{ & m \ar[dl]_{\alpha_m} &   \ar[l]_{\alpha_{m-1}} \\ 1 \ar@(dl, ul)^{\beta} \ar[dr]_{\alpha_1} \\ & 2 \ar[r]_{\alpha_2} & 3 \ar @{{}{*}} @/_2pc/ [uu]}$$
and relations (i) $\alpha_m \cdots \alpha_1 = \beta^2$; (ii) $\overbrace{\alpha_i \cdots \alpha_{i+1} \alpha_i}^{m+1} = 0$ for all $i \in \{1, \ldots, m\} = \mathbb{Z}/\gen{m}$; and (iii) $\alpha_1 \alpha_m = \alpha_1 \beta \alpha_m$ \cite{Asa2}.  Henceforth, we fix $m$, denote this algebra as $A$, and assume $\chr (k) = 2$.  It is well-known that $A$ is simply connected, meaning that it admits no proper connected Galois covers.  In fact, it has no nontrivial radical gradings.  However, there does exist a non-radical $\mathbb{Z}/\gen{2}$-grading of $A$, and we will show that the corresponding smash product is Morita equivalent to a Brauer tree algebra.

Note that $A$ is generated as a $k$-algebra by $\{e_1 + \beta, e_i,\alpha_i\ | \ 1 \leq i \leq  m\}$.  Since relation (iii) can be expressed $\alpha_1 (e_1 + \beta) \alpha_m = 0$ and relation (ii) can be rewritten $\alpha_m \cdots \alpha_1 + (e_1 + \beta)^2 + e_1 = 0$, we see that we obtain a $\mathbb{Z}/\gen{2}$-grading on $A$ with $e_i, \alpha_i \in A_0$ for all $i$ and $e_1 + \beta \in A_1$.  One easily checks that with respect to this grading, $J_G = (\alpha_1, \ldots, \alpha_m)$.  Consequently, the graded simples, up to isomorphism, are the simples $S_i$ concentrated in degree $0$, for $2 \leq i \leq m$, and their shifts $S_i[1]$, as well as the module $e_1 A/ (\alpha_m A + \beta \alpha_m A)$, which is isomorphic to its shift as a graded module.  Likewise, up to isomorphism the indecomposable graded projectives are $e_iA, e_iA[1]$ for $ 2 \leq i \leq m$ and $e_1A \cong e_1A[1]$.

We let $B = A \# k[G]^*$.  From \cite{CoMo}, we know that $J(B) = J_G(A) B$, and hence a $k$-basis of $B/J(B)$ is given by the residue classes of $\{\beta p_g, e_i p_g \mid 1 \leq i \leq m,\  g \in \mathbb{Z}/\gen{2} \}$.  From this, one easily establishes a ring isomorphism $B/J(B) \cong M_2(k) \times k^{2m-2}$ by mapping $$b_0 \beta p_0 + b_1 \beta p_1 + \sum_{i, g} c^g_i e_i p_g \mapsto \left( \left( \begin{array}{cc} c^0_1 & b_1+c^1_1 \\ b_0 + c^0_1 & c^1_1\end{array} \right), c^0_2, c^1_2, \ldots, c^1_m \right).$$
  To compute a basic version of $B$, we can thus take the corner ring $B'$ associated to the full idempotent $1-e_1p_1$.  We make the following observations.
\begin{itemize}
\item $e_i p_1 B e_j p_0 = e_i A_1 e_j p_0 = 0$ for all $i, j > 1$.  
\item $e_i p_0 B e_j p_1 = e_i A_1 e_j p_1 = 0$ for all $i, j > 1$. 
\item $e_i p_g J(B) e_j p_g = e_i J_G(A)_0 e_j p_g$ is $1$-dimensional for all $i, j > 1$ and $g \in \mathbb{Z}/\gen{2}$.
\item $J(B')/J(B')^2$ has a $k$-basis consisting of the residue classes of $\{\alpha_1(e_1 + \beta)p_0, (e_1+\beta)\alpha_m p_1, \alpha_i p_g\ |\ (i, g) \neq (1,1), (m,1) \}$.
\end{itemize}
Hence we can define an isomorphism between $B'$ and the algebra $B''$ given by the quiver 
$$
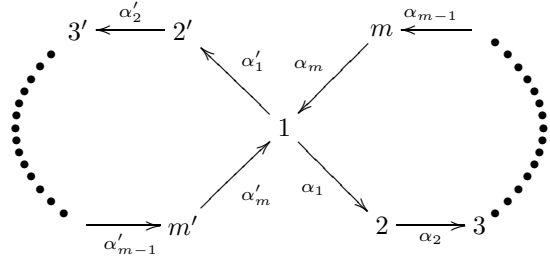
\xymatrix{ 3' \ar@{{}{*}} @/_2pc/ [dd] & 2' \ar[l]_{\alpha'_2} & & m \ar[dl]_{\alpha_m} &   \ar[l]_{\alpha_{m-1}} \\ & & 1 \ar[ul]_{\alpha'_1} \ar[dr]_{\alpha_1} \\ \ar[r]_{\alpha'_{m-1}} & m' \ar[ur]_{\alpha'_m} & & 2 \ar[r]_{\alpha_2} & 3 \ar @{{}{*}} @/_2pc/ [uu]}$$
and relations $\overbrace{\alpha_i \cdots \alpha_{i+1} \alpha_i}^{m+1} = \overbrace{\alpha'_i \cdots \alpha'_{i+1} \alpha'_i}^{m+1}=0$ for all $i$, $\alpha'_1 \alpha_m =  \alpha_1 \alpha'_m = 0$ and $\alpha'_m \cdots \alpha'_1 = \alpha_m \cdots \alpha_1$.  An isomorphism $\varphi : B'' \rightarrow B'$ is given by 
$$\begin{array}{ccc} \begin{array}{cc} \varphi(e_i) = e_i p_0, & 1 \leq i \leq m; \\ \varphi(e_{i'}) = e_i p_1, & 2 \leq i \leq m; \\ \varphi(\alpha_i) = \alpha_i p_0, &  1 \leq i \leq m;  \end{array} & \hspace{5mm} & \varphi(\alpha'_i) = \left\{ \begin{array}{cc} \alpha_i p_1, & 2 \leq i \leq m-1, \\ \alpha_1(e_1 + \beta)p_0, & i = 1, \\ (e_1 + \beta)\alpha_m p_1, & i = m. \end{array} \right. \end{array}$$

\vspace{3mm}
\begin{therm} The nonstandard indecomposable self-injective algebra $A$ of finite representation type and type $(\mathbb{D}_{3m}, 1/3, 1)$ is periodic.  Moreover, its period $p_A$ satisfies $(2m-1) \mid p_A \mid 4(2m-1)$.
\end{therm}

\noindent
{\it Proof.}  We consider $A$ with the $\mathbb{Z}/\gen{2}$-grading described above, and let $B = A \# k[G]^*$, which we have shown is Morita equivalent to a Brauer tree algebra $B'$.  Since $B'$ has $2m-1$ simples and exceptional multiplicity $1$, it is derived equivalent to the symmetric Nakayama algebra with $2m-1$ simples and Loewy length $2m$ \cite{DCSE}, and hence has type $(\mathbb{A}_{2m-1}, \frac{2m-1}{2m-1}, 1)$.  From Table 5.2, we see that the period of $B$ is $ r := 2(2m-1)$.
 
We now follow the strategy of the proof of Theorem 3.7, using special arguments for various details.  First, since $A$ is socle equivalent to the standard algebra of type $(\mathbb{D}_{3m}, 1/3, 1)$, it follows that the the two algebras have the same nonprojective indecomposable modules, and the action of the syzygy functor on objects is the same over either algebra.  Since the standard algebra has period dividing $2(2m-1)$ in this case, all simple $A$-modules are fixed by $\Omega^{r}$.  Hence $\Omega^{r}_{A^e}(A) \cong {}_1 A_{\sigma}$ for some automorphism $\sigma$.  As before, we want to know that $\sigma$ is degree-preserving and that $\Omega^{r}_{A^e}(A)$ is a graded bimodule generated in degree $0$.  Once we have established these facts, we can apply Lemma 3.6 to conclude that $\sigma^2$ is inner.  This yields the stated upper bound of $4(2m-1)$ for the period of $A$, and the lower bound of $(2m-1)$ follows from $p_B | 2p_A$ which can be proved as in Theorem 3.7.

We first address the existence of a graded projective resolution of the bimodule ${}_A A_A$ (note that the argument in Section 2 was for radical gradings only).  In this case, however, since $G$ is abelian, one easily sees that $\bimod[G]{A}$ is equivalent to $\fmod[G]{A^e}$ where $A^e$ is given the grading $A^e_g = \oplus_{h \in G} A_h \otimes A_{h^{-1}g}$.  Thus it is also equivalent to $\fmod{(A^e \#k[G]^*)}$, which clearly has projective covers.  Furthermore, since $J(A^e \#k[G]^*) = J_G(A^e) \cdot (A^e \#k[G]^*) \subseteq J(A^e) \cdot A^e \#k[G]^*$, a graded projective cover will remain a projective cover in the category of (ungraded) bimodules.

As before, that the automorphism $\sigma$ preserves the grading is a consequence of $\Omega^{r}_{A^e}(A)$ being generated in degree $0$.  To prove the latter, we note that the period of $B$ being $r$ implies that $\Omega^{r}(S) \cong S$ as graded modules for each graded simple $A$-module $S$.  It follows that for $P_r = \oplus_{i=1}^m Ae_i \otimes_k e_i A$, each summand with $2 \leq i \leq m$ is generated in degree $0$, while $Ae_1 \otimes_k e_1 A$ can be generated in degree $0$ or $1$ (since $e_iA$ can be).  Thus we see that $\Omega^r_{A^e}(A)$ can be generated in degree $0$ as desired.  $\Box$ \\

As with the standard algebras of torsion order $t>1$, determining the exact periods of the nonstandard algebras is complicated by the difficulty of detecting whether or not the automorphism $\sigma$ in the $2(2m-1)^{th}$ syzygy is inner.  To give one example, when $m = 2$ the standard algebra of type $(D_6, 1/3, 1)$ has period $3$ in characteristic $2$, while the nonstandard algebra $A$ of the same type has period $6$.  In fact, a computation of the beginning of a minimal projective resolution of $A$ yields $\Omega^{3}_{A^e}(A) \cong {}_1 A_{\sigma}$ where $\sigma$ is the automorphism of order $2$ given by $\sigma(\alpha_1) = \alpha_1 (e_1 + \beta), \sigma(\alpha_2) = (e_1+\beta)\alpha_2$ and $\sigma(\beta) = \beta + \beta^2 + \beta^3$.

By the results of Section 4, it follows that the stable Auslander algebra of $A$ is also periodic.  We expect that it is a Galois cover of a deformation of the preprojective algebra $P(\mathbb{D}_{3m})$.  We hope to explore this connection in greater detail, and also investigate similar uses of smash products over deformed preprojective algebras in a future work.

\end{document}